\documentclass[a4paper,10pt]{article}

\usepackage{amsfonts,enumerate,array,amsmath,xcolor,amssymb,latexsym, amsthm, mathrsfs}
\usepackage{mathtools}
\usepackage{booktabs} 
\usepackage{caption} 
\usepackage{subcaption} 
\usepackage{graphicx}
\usepackage[all]{nowidow}
\usepackage[utf8]{inputenc}
\usepackage{tikz}
\usetikzlibrary{er,positioning,bayesnet}
\usepackage{multicol}
\usepackage{algpseudocode,algorithm,algorithmicx}
\usepackage{minted}
\usepackage{hyperref}
\usepackage{bbm}
\usepackage[inline]{enumitem} 
\usepackage{thmtools}
\usepackage{verbatim}
\usepackage[margin=3.0cm]{geometry}
\usepackage{orcidlink}

\numberwithin{equation}{section}

\newtheorem{theorem}{Theorem}[section]
\newtheorem{proposition}[theorem]{Proposition}

\newtheorem{corollary}[theorem]{Corollary}

\newtheorem{rem}[theorem]{Remark}
\theoremstyle{definition}

\newtheorem{hyp}[theorem]{Hypothesis}

\newcommand{\propref}[1]{Proposition~\ref{#1}}

\newcommand{\remref}[1]{Remark~\ref{#1}}

\renewcommand{\eqref}[1]{(\ref{#1})}
\newcommand{\hypref}[1]{Hypothesis~\ref{#1}}
\newcommand{\hyprefs}[1]{Hypotheses~\ref{#1}}

\newcounter{step}

\newcommand{\step}[1]{%
  \refstepcounter{step}
  \textbf{Step \thestep. #1}%
}

\newcommand{\stepref}[1]{%
  \textbf{Step~\ref{#1}}%
}

\allowdisplaybreaks

\usepackage[
backend=biber,
style=numeric,
sorting=nyvt,
maxcitenames=99,
maxbibnames=999
]{biblatex}

\providecommand{\keywords}[1]
{
  \small	
  \noindent\textbf{Keywords.} #1
}

\begin{document}
	
\title{Modelling and analysis of a cocktail phage therapy with coinfection}
\author{Xavier Bardina\orcidlink{0000-0003-1803-3401}$^+$, Marc Cano\orcidlink{0009-0005-6568-3513}$^+$, Silvia Cuadrado\orcidlink{0000-0003-2051-2030}$^+ \ddagger$\footnote{Corresponding author.}
\\{\footnotesize $^+$ Dept. Matem\`atiques, Edifici C, Universitat Aut\`onoma de Barcelona, 08193-Bellaterra}
\\{\footnotesize $^{\ddagger}$ Centre de Recerca Matem\`atica, Edifici C, Campus de Bellaterra, 08193-Bellaterra}
\\{\footnotesize E-mail: Xavier.Bardina@uab.cat, Marc.Cano@uab.cat, Silvia.Cuadrado@uab.cat}
\footnote{This work was supported by the Spanish Ministerio de Ciencia, Innovación y Universidades - Agencia Estatal de Investigación (MICIU/AEI/10.13039/501100011033) under grant No. PID2025-173091NB-I00. M. C. is also supported by the predoctoral grant PRE2022-103868 funded by MCIN/AEI/ 10.13039/501100011033, S. C. has also been supported by the Spanish State Research Agency, through the Severo Ochoa and María de Maeztu Program for Centers and Units of Excellence in R$\&$D (CEX2020-001084-M).}}
\date{\today}
\maketitle

\begin{abstract}
	We analyze a delay differential model for phage therapy in which a cocktail of $n$ different types of bacteriophages, each with its own latency period, is administered and coinfection is allowed. We prove the global existence, uniqueness and nonnegativity of solutions, and the existence of an invariant region, under suitable conditions, in which the bacteria are eradicated exponentially fast. We determine the stability of the bacteria-free equilibrium in terms of the administered doses and give sufficient conditions for the nonexistence of coexistence equilibria, showing that coinfection may allow bacteria and phages to coexist even when the doses are large. For a cocktail of two phage types with the same latency period, we obtain a complete description of the biologically feasible coexistence equilibria, which may be zero, one or two. Numerical simulations illustrate the theoretical results.
\end{abstract}
\keywords{Bacteriophage, Phage therapy, Phage cocktail, Coinfection, Delay differential equations, Stability}


\section{Introduction}

The spread of antibiotic resistance is an important public health concern. As bacterial pathogens acquire resistance mechanisms, an increasing number of infections can no longer be treated with conventional antibiotics. This has renewed the interest in phage therapy, that is, the use of bacteriophages (viruses that infect bacteria) to treat bacterial infections.

Lytic bacteriophages attach to the surface of a bacterial cell and inject their genome into it. Using the machinery of the host cell, they replicate their genetic material and assemble new phage particles. At the end of a latency period, the infected cell lyses and releases the newly formed phages, whose number is known as the burst size, into the surrounding medium, where they can infect other bacteria. Bacteriophages are usually highly specific with respect to the bacteria they can infect. On the one hand, this allows the targeted elimination of a pathogen while minimizing the damage to other bacterial populations. On the other hand, a single phage type may only be effective against a narrow range of bacterial strains. For this reason, phages can be administered as cocktails, that is, mixtures of several phage types, and the design of such cocktails has become an active field of research, see for instance \cite{li2020optimizing,marchi2025multi,van2026rapid,yu2024leveraging}.

Mathematical models have been used to understand the interaction between bacteria and phages. Since the seminal work of Campbell \cite{campbell1961conditions}, the latency period has often been incorporated into these models by means of delay differential equations, see for instance \cite{beretta2001modeling,BARDINA201399,bardina2019coinfection,beldjenna2026assessing}. Other approaches, such as models structured by infection age \cite{calsina2011optimal}, and further analyses of phage-bacteria dynamics can be found in \cite{kyaw2021stability,gomez2023mathematical,wang2024nontrivial}.

Classical single-phage models have been extended in several directions. Models including several phage types have been studied in \cite{bascompte2014models,li2020optimizing,marchi2025multi,van2026rapid,yu2024leveraging}. In particular, in \cite[Chapter~2]{bascompte2014models} a cocktail of $n$ phage types without coinfection is considered, and it is shown that the bacteria-free equilibrium is the only equilibrium when the doses are large enough, whereas for small doses a unique coexistence equilibrium appears. Another extension consists in allowing coinfection \cite{alizon2013co,bardina2019coinfection,sudweeks2025impact}. As in \cite{bardina2019coinfection}, by coinfection here we mean that phages can also adsorb to bacteria that are already infected. Phages adsorbed in this way are lost without producing progeny, which reduces the efficacy of the treatment and may slow down, or even prevent, the eradication of the bacterial population. For a single phage type, it is shown in \cite{bardina2019coinfection} that coinfection substantially changes the equilibria of the system: depending on the burst size and the latency period, a coexistence equilibrium exists either when the dose is large or when it is small.

In this article we consider a model that combines both extensions, namely a cocktail of $n$ phage types in the presence of coinfection, thus generalizing the single-phage model of \cite{bardina2019coinfection} and the coinfection-free cocktail model of \cite{bascompte2014models}. Each phage type has its own adsorption coefficients, burst size, dose and latency period, so that the resulting system involves $n$ different delays. Moreover, coinfection couples the different phage types, since phages of each type can be lost by adsorbing to bacteria infected by phages of any type in the cocktail.

For this model we prove that the initial value problem has a unique global nonnegative solution and that, under suitable conditions on the parameters, the system admits an invariant region in which the bacteria are eradicated exponentially fast. 
We show that the bacteria-free equilibrium always exists and that it is locally asymptotically stable when a weighted sum of the phage doses exceeds the growth rate of the bacteria, and unstable when it is smaller; in particular, the invariant region mentioned above can only exist in the first case.
 We then give sufficient conditions for the nonexistence of coexistence equilibria, which can be interpreted in terms of the expected number of new phages produced by each phage that adsorbs to a bacterium. In particular, we show that coinfection reduces this number, so that effective phages may behave as ineffective ones and, in contrast with the coinfection-free case \cite{bascompte2014models}, coexistence equilibria may exist even when the administered doses are large. We also discuss how combining phage types with different efficacy in the same cocktail affects the existence of such equilibria. In order to obtain existence results and explicit expressions, we study in detail a cocktail of two phage types with the same latency period, for which we obtain a complete description of the biologically feasible coexistence equilibria: depending on the parameters, there may be zero, one or two of them. Finally, numerical simulations illustrate these results and show that the outcome of the therapy may depend on the initial bacterial load: an unstable coexistence equilibrium may act as a threshold between the eradication of the bacteria and their unbounded growth, and the system may be bistable, with small infections being eradicated and larger ones persisting at a positive level.

 The paper is organized as follows. In Section~\ref{The Model} we present the model. Section~\ref{Positivity an Boundedness} is devoted to the existence, uniqueness and nonnegativity of solutions, and in Section~\ref{Invariant Region} we prove the existence of an invariant region. In Section~\ref{Stable States} we study the equilibria of the system, namely the stability of the bacteria-free equilibrium and the existence of coexistence equilibria, and in Subsection~\ref{Particular Case} we focus on the case of two phage types with the same latency period. Numerical simulations are presented in Section~\ref{Numerical Simulations}, and the last section contains some concluding remarks.


\section{The Model}\label{The Model}

In this section we present a model for phage therapy in which a cocktail of $n$ different types of bacteriophages is administered in order to fight a bacterial infection. Let $S(t)$ and $I(t)$ denote the concentrations of susceptible (uninfected) and infected bacteria, respectively, at time $t$, and let $Q_i(t)$, for $i=1,\ldots,n$, denote the concentration of free bacteriophages of type $i$ at time $t$. 

In the absence of phages, susceptible bacteria grow exponentially with net growth rate $\alpha>0$. Phages of type $i$ adsorb to susceptible bacteria according to the law of mass action, with adsorption coefficient $k_i^1>0$, and each adsorption to a susceptible bacterium results in its infection. Hence, susceptible bacteria are infected at rate $\sum_{i=1}^n k_i^1 Q_i(t)S(t)$.

Infected bacteria do not reproduce and die at a natural death rate $\mu>0$. A bacterium infected by a phage of type $i$ lyses after a latency period $\tau_i>0$, provided it has not died before, which happens with probability $e^{-\mu\tau_i}$. Therefore, the rate at which bacteria infected by phages of type $i$ lyse at time $t$ is
\[
k_i^1 e^{-\mu\tau_i} Q_i(t-\tau_i)S(t-\tau_i).
\]

Phages of type $i$ are supplied at a constant rate $d_i>0$, which we refer to as the dose, and decay at rate $m_i>0$. Each lysis of a bacterium infected by a phage of type $i$ releases $b_i>0$ new phages of type $i$, where $b_i$ is the burst size. Phages that adsorb to susceptible bacteria are removed from the free phage population at rate $k_i^1Q_i(t)S(t)$.

Since phages cannot distinguish between infected and uninfected bacteria, phages of type $i$ also adsorb to infected bacteria, with adsorption coefficient $k_i^2>0$, regardless of the type of phage that infected them. We assume that the fate of an infected bacterium is determined by the first phage that infected it: phages that adsorb to an already infected bacterium are lost and modify neither its latency period nor its burst size. Hence, phages of type $i$ are lost at rate $k_i^2Q_i(t)I(t)$, while $I$ is not affected by these adsorptions.

\medskip
These assumptions lead to the following system of delay differential equations:
\begin{equation}\label{system}
	\begin{cases}
		\dot{S}(t)= \displaystyle\Big(\alpha-\sum_{i=1}^n k_i^1Q_i(t)\Big)S(t),\\[2ex]
		\dot{I}(t)= \displaystyle\sum_{i=1}^n k_i^1Q_i(t)S(t)-\mu I(t) - \sum_{i=1}^n k_i^1e^{-\mu\tau_i}Q_i(t-\tau_i)S(t-\tau_i),\\[2ex]
		\dot{Q}_i(t) = d_i-m_iQ_i(t)-k_i^1Q_i(t)S(t)-k_i^2Q_i(t)I(t)+k_i^1b_ie^{-\mu\tau_i}Q_i(t-\tau_i)S(t-\tau_i), \quad i=1,\ldots,n.
	\end{cases}
\end{equation}

Since the delays $\tau_1,\ldots,\tau_n$ may be different, initial data for $S$ must be prescribed on $[-\hat\tau,0]$, where $\hat\tau:=\max_i\tau_i$, and initial data for each $Q_i$ on $[-\tau_i,0]$. No delayed values of $I$ appear in \eqref{system}, so for $I$ it suffices to prescribe $I_0:=I(0)$. However, as is well known for this kind of models (see \cite{bardina2019coinfection,smith2011introduction,smith2012persistence}), this value cannot be chosen arbitrarily. Indeed, the infected bacteria present at time $t$ are exactly those that were infected by a phage of some type $i$ at a time $\theta\in(t-\tau_i,t]$ and have survived until time $t$, that is,
\begin{equation}\label{solI}
	I(t) = \sum_{i=1}^n \int_{t-\tau_i}^{t} k_i^1 e^{-\mu(t-\theta)} Q_i(\theta)S(\theta)\,d\theta .
\end{equation}
For $n=1$, this expression reduces to the one in \cite{smith2011introduction,smith2012persistence,bardina2019coinfection}. Differentiating \eqref{solI}, one checks directly that the right-hand side of \eqref{solI}, which we denote by $J(t)$, satisfies the second equation in \eqref{system}. Conversely, since this equation is linear in $I$, the difference $I(t)-J(t)$ satisfies $\frac{d}{dt}\big(I(t)-J(t)\big)=-\mu\big(I(t)-J(t)\big)$, so that
\[
I(t)=J(t)+\big(I_0-J(0)\big)e^{-\mu t}, \qquad t\geq 0.
\]
Therefore, $I$ is given by \eqref{solI} for all $t\geq 0$ if and only if
\begin{equation}\label{I0}
	I_0 = \sum_{i=1}^n \int_{-\tau_i}^{0} k_i^1 e^{\mu \theta} Q_i^0(\theta) S_0(\theta)\,d\theta
	= \sum_{i=1}^n \int_0^{\tau_i} k_i^1 e^{-\mu s} Q_i^0(-s)S_0(-s)\,ds,
\end{equation}
where $S_0$ and $Q_i^0$ denote the initial data for $S$ and $Q_i$, respectively. We thus impose the following assumption on the initial data.

\begin{hyp}\label{hyp1}
	The initial data are given by continuous nonnegative functions $S_0\colon[-\hat\tau,0]\to\mathbb{R}$ and $Q_i^0\colon[-\tau_i,0]\to\mathbb{R}$, for $i=1,\ldots,n$, and $I_0$ is given by \eqref{I0}.
\end{hyp}


\section{Existence, Uniqueness and Nonnegativity of Solutions}\label{Positivity an Boundedness}

In this section we prove that the initial value problem for \eqref{system} has a unique solution defined for all $t\geq 0$ and that this solution remains nonnegative.

\begin{proposition}\label{prop:existence}
	Under \hypref{hyp1}, the initial value problem for \eqref{system} has a unique solution defined on $[0,\infty)$. Moreover, $S(t)\geq 0$, $I(t)\geq 0$ and $Q_i(t)>0$ for all $t>0$ and $i=1,\ldots,n$.
\end{proposition}

\begin{proof}
	Since the right-hand side of \eqref{system} is locally Lipschitz continuous, standard results for delay differential equations (see \cite{kuang1993delay,smith2011introduction}) ensure the existence of a unique solution defined on a maximal interval $[0,T_{\max})$, with $T_{\max}\in(0,\infty]$. Moreover, if $T_{\max}<\infty$, then the solution is unbounded on $[0,T_{\max})$. We first prove that the solution is nonnegative on $[0,T_{\max})$ and then that $T_{\max}=\infty$.
	
	\medskip
	From the first equation in \eqref{system} we have
	\[
	S(t) = S_0(0) \exp\left( \int_0^t \Big( \alpha - \sum_{i=1}^n k_i^1 Q_i(s) \Big)\, ds \right) \geq 0, \qquad t\in[0,T_{\max}),
	\]
	so that, by \hypref{hyp1}, $S$ is nonnegative on $[-\hat\tau,T_{\max})$.
	
	Now fix $i\in\{1,\ldots,n\}$ and define, for $t\in[0,T_{\max})$,
	\[
	a_i(t) := m_i + k_i^1 S(t) + k_i^2 I(t), \qquad g_i(t) := k_i^1 b_i e^{-\mu\tau_i} Q_i(t-\tau_i) S(t-\tau_i).
	\]
	The equation for $Q_i$ in \eqref{system} reads $\dot Q_i(t) = d_i - a_i(t)Q_i(t) + g_i(t)$, and the variation of constants formula yields
	\begin{equation}\label{voc_Qi}
		Q_i(t) = Q_i^0(0)\, e^{-\int_0^t a_i(r)\,dr} + \int_0^t e^{-\int_s^t a_i(r)\,dr}\,\big(d_i + g_i(s)\big)\, ds .
	\end{equation}
	For $t\in[0,\tau_i]$, $g_i(t)$ only involves the initial data, so $g_i(t)\geq 0$ by \hypref{hyp1}, and \eqref{voc_Qi} together with $d_i>0$ gives $Q_i(t)>0$ for $t\in(0,\tau_i]\cap(0,T_{\max})$. Consequently, $g_i(t)\geq 0$ for $t\in[\tau_i,2\tau_i]\cap[0,T_{\max})$, and \eqref{voc_Qi} gives $Q_i(t)>0$ on this interval as well. Proceeding in this way on the successive intervals $[k\tau_i,(k+1)\tau_i]$, $k\in\mathbb{N}$, we obtain $Q_i(t)>0$ for all $t\in(0,T_{\max})$.
	
	Finally, as shown in Section~\ref{The Model}, \hypref{hyp1} implies that $I$ is given by \eqref{solI}. Since the integrands in \eqref{solI} are nonnegative, we conclude that $I(t)\geq 0$ for all $t\in[0,T_{\max})$.
	
	\medskip
	In order to prove global existence let us assume, by contradiction, that $T:=T_{\max}<\infty$. We show that the solution is bounded on $[0,T)$. Since $Q_i\geq 0$ for all $i$, we have $\dot S(t)\leq \alpha S(t)$, and hence $S(t)\leq S_0(0)e^{\alpha t}$ for $t\in[0,T)$. Therefore,
	\[
	S(t)\leq C_T:=\max\Big\{\max_{s\in[-\hat\tau,0]} S_0(s),\; S_0(0)e^{\alpha T}\Big\}, \qquad t\in[-\hat\tau,T).
	\]
	Using that $S$, $I$ and $Q_i$ are nonnegative, the equation for $Q_i$ gives
	\[
	\dot Q_i(t) \leq d_i + K_i\, Q_i(t-\tau_i), \qquad t\in[0,T),
	\]
	where $K_i:=k_i^1 b_i e^{-\mu\tau_i} C_T$. Integrating on $[0,t]$ and performing the change of variables $u=s-\tau_i$, we obtain
	\[
	Q_i(t) \leq Q_i^0(0) + d_i t + K_i\int_{-\tau_i}^{t-\tau_i} Q_i(u)\,du
	\leq c_i + K_i\int_0^{t} Q_i(u)\,du, \qquad t\in[0,T),
	\]
	where $c_i := Q_i^0(0) + d_i T + K_i\tau_i \max_{s\in[-\tau_i,0]} Q_i^0(s)$. By Gronwall's lemma, $Q_i(t)\leq c_i e^{K_i t}\leq c_i e^{K_i T}$ for $t\in[0,T)$. Finally, from \eqref{solI} and the previous bounds,
	\[
	I(t)\leq C_T\sum_{i=1}^n k_i^1\tau_i \max\Big\{\max_{s\in[-\tau_i,0]}Q_i^0(s),\; c_i e^{K_i T}\Big\}, \qquad t\in[0,T).
	\]
	Thus, the solution is bounded on $[0,T_{\max})$, which contradicts $T_{\max}<\infty$. Hence $T_{\max}=\infty$.
\end{proof}


\section{Invariant Region}\label{Invariant Region}

In this section we show that, under suitable conditions on the parameters and on the initial data, the solutions of \eqref{system} remain in a bounded region. Following \cite{BARDINA201399,bardina2019coinfection}, for each phage type $i$ we fix a constant $M_i>0$, which we call the maximal attack rate of phages of type $i$ and which can be interpreted as the phage concentration above which the infection process becomes saturated. We define
\begin{equation}\label{Sbar}
	\bar S \coloneqq \min_{i=1,\ldots,n} \frac{m_iM_i-d_i}{k_i^1 b_i e^{-\mu\tau_i} M_i},
	\qquad
	\eta \coloneqq \frac{1}{\mu}\sum_{j=1}^n \frac{m_jM_j-d_j}{b_j e^{-\mu\tau_j}},
\end{equation}
and, for $i=1,\ldots,n$,
\begin{equation}\label{gamma_i}
	\gamma_i \coloneqq \frac{d_i}{m_i + k_i^2\eta}.
\end{equation}

\begin{hyp}\label{hyp2}
	For every $i=1,\ldots,n$ we have $d_i<m_iM_i$ and $b_ie^{-\mu\tau_i}>1$. Moreover,
	\[
	\alpha \leq \sum_{i=1}^n k_i^1\gamma_i .
	\]
\end{hyp}

Under \hypref{hyp2}, we have $\bar S>0$, $\eta>0$ and $0<\gamma_i<d_i/m_i<M_i$ for every $i$, so that the region
\[
R \coloneqq [0,\bar S]\times[0,\eta]\times[\gamma_1,M_1]\times\cdots\times[\gamma_n,M_n]
\]
is nonempty.

\begin{hyp}\label{hyp3}
	The initial data satisfy the following conditions:
	\begin{enumerate}[label=(\roman*)]
		\item $S_0(s)<\bar S$ for $s\in[-\hat\tau,0]$, $\gamma_i\leq Q_i^0(s)\leq M_i$ for $s\in[-\tau_i,0]$ and $i=1,\ldots,n$, and $I_0\leq\eta$.
		\item For every $i=1,\ldots,n$ and $s\in[-\tau_i,0]$,
		\[
		b_ie^{-\mu\tau_i}\,Q_i^0(s)S_0(s) > \gamma_i S_0(0).
		\]
	\end{enumerate}
\end{hyp}

Note that \hypref{hyp3}~(ii) implies $S_0(0)>0$. Consequently, since $\gamma_i>0$, the same condition yields $S_0(s)>0$ for all $s\in[-\tau_i,0]$ and every $i$, that is, $S_0(s)>0$ for all $s\in[-\hat\tau,0]$.

\begin{rem}\label{rem:regime}
	Since $\gamma_i<d_i/m_i$ for every $i$, \hypref{hyp2} implies
	\[
	\alpha < \sum_{i=1}^n \frac{k_i^1 d_i}{m_i}.
	\]
	Hence, the invariant region obtained in this section only exists in the regime in which the bacteria-free equilibrium is locally asymptotically stable (see Proposition~\ref{prop:E0}).
\end{rem}

\begin{proposition}\label{prop:invariant}
	Under \hyprefs{hyp1}, \ref{hyp2} and \ref{hyp3}, the solution of \eqref{system} satisfies
	\[
	\big(S(t),I(t),Q_1(t),\ldots,Q_n(t)\big)\in R \qquad\text{for all } t\geq 0 .
	\]
	Moreover, $S$ is nonincreasing on $[0,\infty)$.
\end{proposition}

\begin{proof}
	Let $x(t)\coloneqq\big(S(t),I(t),Q_1(t),\ldots,Q_n(t)\big)$. We first note that, by \eqref{Sbar} and \eqref{gamma_i},
	\begin{equation}\label{identities_R}
		d_i - m_i\gamma_i - k_i^2\gamma_i\eta = 0
		\qquad\text{and}\qquad
		k_i^1 b_i e^{-\mu\tau_i} M_i\, \bar S \leq m_iM_i-d_i, \qquad i=1,\ldots,n,
	\end{equation}
	and that $\sum_{i=1}^n (m_iM_i-d_i)/(b_ie^{-\mu\tau_i}) = \mu\eta$. Recall also that, by \propref{prop:existence}, the solution exists for all $t\geq 0$, and $S$ and $I$ are nonnegative.
	
	Consider the set
	\[
	A\coloneqq\big\{T\geq 0 \colon x(t)\in R \text{ for all } t\in[0,T]\big\}.
	\]
	By \hypref{hyp3}~(i), $x(0)\in R$, so $0\in A$. Let $T^*\coloneqq\sup A$ and assume, by contradiction, that $T^*<\infty$. Since $R$ is closed and $x$ is continuous, $T^*\in A$. We will show that there exists $\delta>0$ such that $T^*+\delta\in A$, which contradicts the definition of $T^*$.
	
	\medskip
	\step{$S$ is nonincreasing on $[0,T^*]$.}\label{step1}
    
	For $t\in[0,T^*]$ we have $Q_i(t)\geq\gamma_i$ for every $i$ and hence, by \hypref{hyp2},
	\[
	\dot S(t) \leq \Big(\alpha-\sum_{i=1}^n k_i^1\gamma_i\Big)S(t)\leq 0 .
	\]
	Therefore, $S$ is nonincreasing on $[0,T^*]$. Since $S_0(0)>0$, the explicit expression of $S$ in the proof of \propref{prop:existence} gives $S(t)>0$, and thus $0<S(t)\leq S_0(0)<\bar S$ for $t\in[0,T^*]$. Together with \hypref{hyp3}~(i), this yields
	\begin{equation}\label{S_strict}
		S(s)<\bar S \qquad\text{for all } s\in[-\hat\tau,T^*].
	\end{equation}
	
	\medskip
	\step{If $Q_i(T^*)=M_i$, then $\dot Q_i(T^*)<0$.}\label{step2}
    
	We have $Q_i(T^*-\tau_i)\leq M_i$, by \hypref{hyp3}~(i) if $T^*<\tau_i$ and because $T^*\in A$ otherwise. Dropping the nonpositive terms $-k_i^1Q_iS$ and $-k_i^2Q_iI$ in the equation for $Q_i$ and using \eqref{S_strict} and \eqref{identities_R}, we obtain
	\[
	\dot Q_i(T^*) \leq d_i - m_iM_i + k_i^1 b_i e^{-\mu\tau_i} M_i\, S(T^*-\tau_i)
	< d_i - m_iM_i + k_i^1 b_i e^{-\mu\tau_i} M_i\, \bar S \leq 0 .
	\]
	
	\medskip
	\step{If $Q_i(T^*)=\gamma_i$, then $\dot Q_i(T^*)>0$.}\label{step3}
    
	Since $I(T^*)\leq\eta$, the first identity in \eqref{identities_R} gives
	\begin{eqnarray*}
	&&\dot Q_i(T^*) \geq d_i - m_i\gamma_i - k_i^2\gamma_i\eta - k_i^1\gamma_i S(T^*) + k_i^1 b_i e^{-\mu\tau_i} Q_i(T^*-\tau_i)S(T^*-\tau_i)\\
	&=& k_i^1\Big( b_i e^{-\mu\tau_i} Q_i(T^*-\tau_i)S(T^*-\tau_i) - \gamma_i S(T^*)\Big).
	\end{eqnarray*}
	If $T^*<\tau_i$, then $T^*-\tau_i\in[-\tau_i,0)$ and $S(T^*)\leq S_0(0)$, so that, by \hypref{hyp3}~(ii),
	\[
	\dot Q_i(T^*) \geq k_i^1\Big( b_i e^{-\mu\tau_i} Q_i^0(T^*-\tau_i)S_0(T^*-\tau_i) - \gamma_i S_0(0)\Big) > 0 .
	\]
	If $T^*\geq\tau_i$, then $T^*-\tau_i\in[0,T^*]$, so that $Q_i(T^*-\tau_i)\geq\gamma_i$ and, by \stepref{step1}, $S(T^*-\tau_i)\geq S(T^*)>0$. Therefore, by \hypref{hyp2},
	\[
	\dot Q_i(T^*) \geq k_i^1\big(b_ie^{-\mu\tau_i}-1\big)\gamma_i S(T^*) > 0 .
	\]
	
	\medskip
	\step{If $I(T^*)=\eta$, then $\dot I(T^*)<0$.}\label{step4}
    
	Since $Q_i(T^*)\leq M_i$ and $S(T^*)<\bar S$, the second inequality in \eqref{identities_R} yields
	\[
	\sum_{i=1}^n k_i^1Q_i(T^*)S(T^*) < \sum_{i=1}^n k_i^1 M_i\,\bar S \leq \sum_{i=1}^n \frac{m_iM_i-d_i}{b_ie^{-\mu\tau_i}} = \mu\eta .
	\]
	Hence,
	\[
	\dot I(T^*) < \mu\eta - \mu\eta - \sum_{i=1}^n k_i^1 e^{-\mu\tau_i}Q_i(T^*-\tau_i)S(T^*-\tau_i) \leq 0 .
	\]
	
	\medskip

To conclude, by continuity and Step 2, Step 3 and Step 4, there exists $\delta>0$ such that $\gamma_i\leq Q_i(t)\leq M_i$ and $I(t)\leq\eta$ for all $t\in[0,T^*+\delta]$ and every $i$: the inequalities that hold strictly at $T^*$ persist by continuity, and those that hold with equality persist because of the strict sign of the corresponding derivative. Arguing as in Step 1, we then have $0\leq S(t)\leq S_0(0)<\bar S$ for $t\in[0,T^*+\delta]$, so that $T^*+\delta\in A$, which contradicts the definition of $T^*$. Hence $T^*=\infty$ and, by Step 1, $S$ is nonincreasing on $[0,\infty)$.
\end{proof}

If the inequality in \hypref{hyp2} is strict, the previous result shows that the bacteria are eradicated exponentially fast.

\begin{corollary}\label{cor:eradication}
	Assume \hyprefs{hyp1}, \ref{hyp2} and \ref{hyp3}, and suppose moreover that $c\coloneqq\sum_{i=1}^n k_i^1\gamma_i-\alpha>0$. Then
	\[
	S(t)\leq S_0(0)e^{-ct} \quad\text{for all } t\geq 0,
	\qquad\text{and}\qquad
	I(t)\leq S_0(0)\Big(\sum_{i=1}^n k_i^1M_i\tau_i\Big)e^{-c(t-\hat\tau)} \quad\text{for all } t\geq\hat\tau .
	\]
	In particular, $S(t)\to 0$ and $I(t)\to 0$ as $t\to\infty$.
\end{corollary}

\begin{proof}
	By \propref{prop:invariant}, $Q_i(t)\geq\gamma_i$ for all $t\geq0$ and every $i$, so that $\dot S(t)\leq -cS(t)$, which gives the bound for $S$. For $t\geq\hat\tau$, using \eqref{solI}, $Q_i\leq M_i$ and the fact that $S$ is nonincreasing on $[0,\infty)$, we obtain
	\[
	I(t)\leq \sum_{i=1}^n k_i^1 M_i \int_{t-\tau_i}^{t} S(\theta)\,d\theta \leq \Big(\sum_{i=1}^n k_i^1M_i\tau_i\Big) S(t-\hat\tau) \leq S_0(0)\Big(\sum_{i=1}^n k_i^1M_i\tau_i\Big)e^{-c(t-\hat\tau)} .
	\]
\end{proof}


\section{Equilibria and Stability}\label{Stable States}

This section is devoted to the study of the equilibria of system \eqref{system}, that is, the solutions of
\begin{equation}\label{stability}
	\left\lbrace
	\begin{aligned}
		0 &= \Big(\alpha-\sum_{i=1}^nk_i^1Q_i\Big)S,\\
		0 &= \sum_{i=1}^nk_i^1Q_iS-\mu I-\sum_{i=1}^nk_i^1e^{-\mu\tau_i}Q_iS,\\
		0 &= d_i - m_iQ_i-k_i^1Q_iS - k_i^2Q_iI+k_i^1b_ie^{-\mu\tau_i}Q_iS, \qquad i=1,\ldots,n.
	\end{aligned}
	\right.
\end{equation}
From the first equation, either $S=0$ or
\begin{equation}\label{equation_alpha}
	\alpha = \sum_{i=1}^nk_i^1 Q_i.
\end{equation}
If $S=0$, the remaining equations give $I=0$ and $Q_i=d_i/m_i$ for $i=1,\ldots,n$, which yields the bacteria-free equilibrium
\[
E_0 = \Big(0,0,\frac{d_1}{m_1},\ldots,\frac{d_n}{m_n}\Big),
\]
which exists for any choice of the parameters. Equilibria with $S>0$ are called coexistence equilibria, and we say that a coexistence equilibrium is biologically feasible if all its components are positive.

\begin{proposition}\label{prop:E0}
	The bacteria-free equilibrium $E_0$ is locally asymptotically stable if
	\[
	\alpha < \sum_{i=1}^n\frac{k_i^1d_i}{m_i}
	\]
	and unstable if the reverse strict inequality holds.
\end{proposition}

\begin{proof}
	We linearize \eqref{system} around $E_0$. Writing $S(t)=u(t)$, $I(t)=v(t)$ and $Q_i(t)=d_i/m_i+w_i(t)$, $i=1,\ldots,n$, and neglecting the terms of order higher than one, we obtain
	\begin{align*}
		\dot u(t) &= \lambda_0\, u(t), \qquad\text{where } \lambda_0\coloneqq\alpha-\sum_{i=1}^n\frac{k_i^1d_i}{m_i},\\
		\dot v(t) &= \sum_{i=1}^n \frac{k_i^1d_i}{m_i}\Big(u(t)-e^{-\mu\tau_i}u(t-\tau_i)\Big) - \mu v(t),\\
		\dot w_i(t) &= -m_iw_i(t) - \frac{k_i^1d_i}{m_i}\Big(u(t)-b_ie^{-\mu\tau_i}u(t-\tau_i)\Big) - \frac{k_i^2d_i}{m_i}\,v(t), \qquad i=1,\ldots,n.
	\end{align*}
	Looking for solutions of the form $e^{\lambda t}(u,v,w_1,\ldots,w_n)$ with $\lambda\in\mathbb{C}$, we obtain a linear system which is lower triangular, with diagonal entries $\lambda-\lambda_0$, $\lambda+\mu$, $\lambda+m_1,\ldots,\lambda+m_n$, and the characteristic equation is
	\[
	(\lambda-\lambda_0)(\lambda+\mu)\prod_{i=1}^n(\lambda+m_i) = 0 .
	\]
	Its roots are $\lambda_0$, $-\mu$, $-m_1,\ldots,-m_n$. Since all of them except possibly $\lambda_0$ are negative, standard results for delay differential equations (see \cite{kuang1993delay,smith2011introduction}) imply that $E_0$ is locally asymptotically stable if $\lambda_0<0$ and unstable if $\lambda_0>0$.
\end{proof}

In view of \propref{prop:E0}, we will say that the doses are large if $\sum_{i=1}^n k_i^1d_i/m_i>\alpha$ and small if $\sum_{i=1}^n k_i^1d_i/m_i<\alpha$. Thus, the bacteria-free equilibrium is locally asymptotically stable when the doses are large and unstable when they are small.

Let us now consider a coexistence equilibrium $(S,I,Q_1,\ldots,Q_n)$, so that $S>0$ and \eqref{equation_alpha} holds. From the second equation in \eqref{stability} we obtain
\begin{equation}\label{equation_I}
	I = B S, \qquad\text{where } B\coloneqq \frac{1}{\mu}\sum_{j=1}^n k_j^1 Q_j\big(1-e^{-\mu\tau_j}\big).
\end{equation}
Substituting \eqref{equation_I} into the equation for $Q_i$ we get $d_i = Q_i\,(m_i + B_iS)$, where
\begin{equation}\label{B_i}
	B_i \coloneqq -k_i^1\big(b_ie^{-\mu\tau_i}-1\big) + k_i^2 B, \qquad i=1,\ldots,n.
\end{equation}
Since $d_i>0$, we have $m_i+B_iS\neq 0$,
\begin{equation}\label{equation_Qi}
	Q_i = \frac{d_i}{m_i + B_iS},
\end{equation}
and $Q_i>0$ if and only if $m_i+B_iS>0$. Substituting \eqref{equation_Qi} into \eqref{equation_alpha}, we obtain
\begin{equation}\label{F}
	\alpha = \sum_{i=1}^n \frac{k_i^1 d_i}{m_i + B_iS} \eqqcolon F(S,Q_1,\ldots,Q_n),
\end{equation}
where the dependence on $Q_1,\ldots,Q_n$ comes from \eqref{equation_I}. Notice that $F(0,Q_1,\ldots,Q_n)=\sum_{i=1}^n k_i^1d_i/m_i$. Moreover, if the coefficients $B_1,\ldots,B_n$ are regarded as fixed, then
\begin{equation}\label{derivative_F}
	\frac{\partial F}{\partial S} = -\sum_{i=1}^n \frac{k_i^1 d_i B_i}{(m_i + B_iS)^2},
\end{equation}
so that $F$ is decreasing in $S$ if $B_i>0$ for every $i$, and increasing in $S$ on the range where $m_i+B_iS>0$ for every $i$ if $B_i<0$ for every $i$. In the first case, equation \eqref{F} can only have a positive solution if $F(0,Q_1,\ldots,Q_n)>\alpha$, and in the second case only if $F(0,Q_1,\ldots,Q_n)<\alpha$. However, the coefficients $B_i$ depend on the equilibrium through $B$, and therefore their signs are not known a priori. This leads us to look for bounds of $B_i$ that only depend on the parameters of the model.

Let $\underline\tau\coloneqq\min_i\tau_i$ and recall that $\hat\tau=\max_i\tau_i$. If the coexistence equilibrium is biologically feasible, then $Q_j>0$ for every $j$ and, by \eqref{equation_alpha},
\[
\frac{\alpha}{\mu}\big(1-e^{-\mu\underline\tau}\big) \leq B \leq \frac{\alpha}{\mu}\big(1-e^{-\mu\hat\tau}\big).
\]
Consequently, $\underline A_i \leq B_i \leq \overline A_i$ for $i=1,\ldots,n$, where
\begin{equation}\label{A_bounds}
	\underline A_i \coloneqq -k_i^1\big(b_ie^{-\mu\tau_i}-1\big) + \frac{k_i^2\alpha}{\mu}\big(1-e^{-\mu\underline\tau}\big),
	\qquad
	\overline A_i \coloneqq -k_i^1\big(b_ie^{-\mu\tau_i}-1\big) + \frac{k_i^2\alpha}{\mu}\big(1-e^{-\mu\hat\tau}\big).
\end{equation}
Note that these bounds only depend on the parameters of the model.

\begin{proposition}\label{Statement1}
	\begin{enumerate}[label=(\roman*)]
		\item Assume that $\underline A_i>0$ for every $i=1,\ldots,n$. If there exists a biologically feasible coexistence equilibrium, then $\alpha<\sum_{i=1}^n k_i^1d_i/m_i$. Equivalently, if $\sum_{i=1}^n k_i^1d_i/m_i\leq\alpha$, there is no biologically feasible coexistence equilibrium.
		\item Assume that $\overline A_i<0$ for every $i=1,\ldots,n$. If there exists a biologically feasible coexistence equilibrium, then $\alpha>\sum_{i=1}^n k_i^1d_i/m_i$. Equivalently, if $\sum_{i=1}^n k_i^1d_i/m_i\geq\alpha$, there is no biologically feasible coexistence equilibrium.
	\end{enumerate}
\end{proposition}

\begin{proof}
	Let $(S,I,Q_1,\ldots,Q_n)$ be a biologically feasible coexistence equilibrium, so that $S>0$ and $m_i+B_iS>0$ for every $i$. In case (i), $B_i\geq\underline A_i>0$ for every $i$, so that, by \eqref{derivative_F}, $F$ is strictly decreasing in $S$ and
	\[
	\alpha = F(S,Q_1,\ldots,Q_n) < F(0,Q_1,\ldots,Q_n) = \sum_{i=1}^n \frac{k_i^1 d_i}{m_i}.
	\]
	In case (ii), $B_i\leq\overline A_i<0$ for every $i$, so that $F$ is strictly increasing in $S$ on $[0,S]$, where $m_i+B_iS>0$, and the reverse strict inequality follows.
\end{proof}

We say that a phage of type $i$ is effective if $b_ie^{-\mu\tau_i}>1$, that is, if the expected number of phages released by a bacterium infected by a phage of type $i$, taking into account the probability that it survives the latency period, is larger than one. If $b_i>1$, this is equivalent to $\tau_i<\tau_i^*\coloneqq\frac{1}{\mu}\log b_i$, while phages with $b_i\leq 1$ would never be effective. Since the second term in $\underline A_i$ is positive, the condition $\underline A_i>0$ holds automatically for ineffective phages. Therefore, the hypothesis of \propref{Statement1}~(i) only has to be checked for the effective phages of the cocktail; for them, it holds when the latency period $\tau_i$ is close enough to $\tau_i^*$ or the coinfection coefficient $k_i^2$ is large enough. On the other hand, the hypothesis of \propref{Statement1}~(ii) requires all the phages of the cocktail to be effective.

\begin{rem}\label{rem:R0}
	The conditions in \propref{Statement1} admit the following interpretation: for a biologically feasible coexistence equilibrium, let
	\[
	R_{0,i} \coloneqq \frac{k_i^1b_ie^{-\mu\tau_i}}{k_i^1+k_i^2B}, \qquad i=1,\ldots,n.
	\]
	Since $I=BS$, the quantity $k_i^1/(k_i^1+k_i^2B) = k_i^1S/(k_i^1S+k_i^2I)$ is the probability that a phage of type $i$ adsorbing to a bacterium adsorbs to a susceptible one, while $b_ie^{-\mu\tau_i}$ is the expected number of phages of type $i$ released by a bacterium infected by a phage of type $i$. Hence, $R_{0,i}$ is the expected number of phages of type $i$ produced per adsorbed phage of type $i$, without taking into account the decay of phages. Since $B_i=(k_i^1+k_i^2B)(1-R_{0,i})$, the sign of $B_i$ is the sign of $1-R_{0,i}$, and the bounds on $B$ yield
	\[
	\frac{k_i^1b_ie^{-\mu\tau_i}}{k_i^1+\frac{k_i^2\alpha}{\mu}(1-e^{-\mu\hat\tau})} \leq R_{0,i} \leq \frac{k_i^1b_ie^{-\mu\tau_i}}{k_i^1+\frac{k_i^2\alpha}{\mu}(1-e^{-\mu\underline\tau})}.
	\]
	The hypothesis of \propref{Statement1}~(i) means that the upper bound is smaller than one for every $i$, that is, no phage type can replace itself through infection, while the hypothesis of \propref{Statement1}~(ii) means that the lower bound is larger than one for every $i$.
\end{rem}

\begin{rem}\label{rem:n1}
	For a single phage type ($n=1$), the previous analysis is complete. In this case $\underline\tau=\hat\tau=\tau_1$, so that $B=\alpha(1-e^{-\mu\tau_1})/\mu$ and $B_1=\underline A_1=\overline A_1\eqqcolon A$ only depends on the parameters. If $A\neq 0$, equation \eqref{F} has the unique solution
	\[
	S^* = \frac{1}{A}\Big(\frac{k_1^1d_1}{\alpha}-m_1\Big),
	\]
	and $m_1+AS^*=k_1^1d_1/\alpha>0$. Hence, there exists a biologically feasible coexistence equilibrium if and only if $A\big(k_1^1d_1/m_1-\alpha\big)>0$, and in this case it is unique. Moreover, $A$ is an increasing function of $\tau_1$. If $b_1>1$, it vanishes at
	\[
	\tilde\tau \coloneqq \frac{1}{\mu}\log\left(\frac{k_1^1b_1+\frac{k_1^2\alpha}{\mu}}{k_1^1+\frac{k_1^2\alpha}{\mu}}\right) < \tau_1^*,
	\]
	while if $b_1\leq 1$ it is positive for every $\tau_1$. Thus, for $\tau_1>\tilde\tau$ the coexistence equilibrium exists if and only if the dose is large, and for $\tau_1<\tilde\tau$ it exists if and only if the dose is small. In the absence of coinfection ($k_1^2=0$), the threshold $\tilde\tau$ coincides with $\tau_1^*$, so that coinfection lowers the latency period above which coexistence occurs for large doses. This is in agreement with the analysis of the single-phage model with coinfection in \cite{bardina2019coinfection}.
\end{rem}

Let us discuss the consequences of these results. Under the hypothesis of \propref{Statement1}~(i), if the doses are small, the bacteria-free equilibrium is unstable and there is no biologically feasible coexistence equilibrium, which suggests that the phage cocktail is not able to control the bacterial population. If the doses are large, the bacteria-free equilibrium is locally asymptotically stable but coexistence equilibria are not excluded; in Section~\ref{Particular Case} we will show that they indeed exist in a particular case, so that the outcome of the therapy may depend on the initial conditions. Under the hypothesis of \propref{Statement1}~(ii), if the doses are large, the bacteria-free equilibrium is locally asymptotically stable and it is the only biologically feasible equilibrium, whereas if the doses are small it is unstable and coexistence equilibria may exist.

These results can be compared with those in \cite[Chapter~2]{bascompte2014models}, where a cocktail of $n$ phage types without coinfection is considered and it is shown that, when the doses are large, the bacteria-free equilibrium is the only equilibrium, while when the doses are small there exists a unique coexistence equilibrium. In the absence of coinfection ($k_i^2=0$ for all $i$), we have $B_i=-k_i^1(b_ie^{-\mu\tau_i}-1)$, so that, if all the phages are effective, we are under the hypothesis of \propref{Statement1}~(ii), in agreement with \cite{bascompte2014models}. Coinfection adds the positive term $k_i^2B$ to $B_i$, that is, it reduces $R_{0,i}$, since phages are lost by adsorbing to infected bacteria. As a consequence, a cocktail of effective phages may satisfy the hypothesis of \propref{Statement1}~(i) and behave as a cocktail of ineffective phages: coexistence equilibria can then exist even when the doses are large.

Finally, since the hypothesis of \propref{Statement1}~(i) holds automatically for ineffective phages, a cocktail combining effective and ineffective phages satisfies it as soon as $\underline A_i>0$ for each effective phage. In this situation, adding phages of any type to the cocktail increases $\sum_{i=1}^n k_i^1d_i/m_i$ and can therefore turn small doses into large ones. This makes the bacteria-free equilibrium locally asymptotically stable but, at the same time, coexistence equilibria are no longer excluded.

When all the latency periods coincide, $\underline A_i=\overline A_i$ and $B_i$ only depends on the parameters of the model. This is the situation considered in the next subsection.


\subsection{Particular case: two phage types with the same latency period}\label{Particular Case}

\propref{Statement1} gives conditions under which coexistence equilibria cannot exist, but it does not guarantee their existence. In order to obtain existence results and explicit expressions for the coexistence equilibria, in this subsection we consider a cocktail of two phage types ($n=2$) with the same latency period $\tilde\tau\coloneqq\tau_1=\tau_2$.

In this case, \eqref{equation_alpha} implies that every coexistence equilibrium satisfies $B=\frac{\alpha}{\mu}(1-e^{-\mu\tilde\tau})$, so that the coefficients $B_i$ in \eqref{B_i} only depend on the parameters of the model. More precisely, $B_i=\underline A_i=\overline A_i=A_i$, where
\begin{equation}\label{A_i}
	A_i\coloneqq -k_i^1\big(b_ie^{-\mu\tilde\tau}-1\big)+\frac{k_i^2\alpha}{\mu}\big(1-e^{-\mu\tilde\tau}\big), \qquad i=1,2.
\end{equation}
We assume $A_1A_2\neq 0$; the cases in which one of the coefficients vanishes can be treated in a similar way. By \eqref{equation_I} and \eqref{equation_Qi}, a coexistence equilibrium is determined by its first component $S>0$ through
\begin{equation}\label{Q_i}
	I=\frac{\alpha}{\mu}\big(1-e^{-\mu\tilde\tau}\big)S, \qquad Q_i=\frac{d_i}{m_i+A_iS}, \quad i=1,2,
\end{equation}
and, by \eqref{F}, $S$ is a solution of the equation $F(S)=\alpha$, where now
\begin{equation}\label{F_particular}
	F(S) = \frac{k_1^1d_1}{m_1+A_1S} + \frac{k_2^1d_2}{m_2+A_2S}
\end{equation}
only depends on $S$. Conversely, every solution $S>0$ of $F(S)=\alpha$ yields a coexistence equilibrium through \eqref{Q_i}. This equilibrium is biologically feasible if and only if $m_i+A_iS>0$ for $i=1,2$, that is, if and only if $S$ belongs to the interval
\[
D\coloneqq\big(0,S_{\max}\big), \qquad S_{\max}\coloneqq\min\Big\{-\frac{m_i}{A_i}\colon A_i<0\Big\},
\]
with the convention $S_{\max}=\infty$ if $A_1,A_2>0$. Thus, biologically feasible coexistence equilibria correspond to solutions in $D$ of $F(S)=\alpha$. On $D$, the function $F$ is smooth and positive, $F(0)=\frac{k_1^1d_1}{m_1}+\frac{k_2^1d_2}{m_2}$, and
\begin{equation}\label{derivative_F_particular}
	F'(S) = -\frac{k_1^1d_1A_1}{(m_1+A_1S)^2}-\frac{k_2^1d_2A_2}{(m_2+A_2S)^2},
	\qquad
	F''(S) = \frac{2k_1^1d_1A_1^2}{(m_1+A_1S)^3}+\frac{2k_2^1d_2A_2^2}{(m_2+A_2S)^3}>0 .
\end{equation}

In order to obtain explicit expressions, note that a direct computation gives
\begin{equation}\label{nova}
(m_1+A_1S)(m_2+A_2S)\big(F(S)-\alpha\big) = -P(S),
\end{equation}
where
\begin{equation}\label{equacio_quadratica_particular}
	P(S)\coloneqq\alpha A_1A_2S^2+\big(A_1(\alpha m_2-k_2^1d_2)+A_2(\alpha m_1-k_1^1d_1)\big)S+m_1m_2\big(\alpha-F(0)\big).
\end{equation}
Hence, the solutions of $F(S)=\alpha$ are roots of the quadratic polynomial $P$, and in particular there are at most two of them. When the discriminant
\begin{equation}\label{L}
	L\coloneqq \big(A_1(\alpha m_2-k_2^1d_2)+A_2(\alpha m_1-k_1^1d_1)\big)^2-4\alpha A_1A_2\,m_1m_2\big(\alpha-F(0)\big)
\end{equation}
is nonnegative, these roots are
\begin{equation}\label{S_hat}
	\hat S\coloneqq\frac{A_1(k_2^1d_2-\alpha m_2)+A_2(k_1^1d_1-\alpha m_1)+\sqrt{L}}{2\alpha A_1A_2}
	\quad\text{and}\quad
	S'\coloneqq\frac{A_1(k_2^1d_2-\alpha m_2)+A_2(k_1^1d_1-\alpha m_1)-\sqrt{L}}{2\alpha A_1A_2}.
\end{equation}
Notice that $\hat S\geq S'$ if $A_1A_2>0$ and $\hat S\leq S'$ if $A_1A_2<0$.

\begin{proposition}\label{prop:particular}
	Let us assume that $A_1$ and $A_2$ defined in \eqref{A_i} satisfy $A_1A_2\neq0$.
	\begin{enumerate}[label=(\roman*)]
		\item If $A_1>0$ and $A_2>0$, there exists a biologically feasible coexistence equilibrium if and only if $F(0)>\alpha$. In this case, it is unique and its first component is $\hat S$.
		\item If $A_1<0$ and $A_2<0$, there exists a biologically feasible coexistence equilibrium if and only if $F(0)<\alpha$. In this case, it is unique and its first component is $S'$.
		\item Assume $A_1A_2<0$.
		\begin{enumerate}[label=(\alph*)]
			\item If $F(0)<\alpha$, there exists a unique biologically feasible coexistence equilibrium, and its first component is $S'$.
			\item If $F(0)>\alpha$ and $F'(0)\geq0$, there is no biologically feasible coexistence equilibrium.
			\item If $F(0)>\alpha$ and $F'(0)<0$, the function $F$ attains its minimum on $D$ at a unique point $S_m$. There is no biologically feasible coexistence equilibrium if $F(S_m)>\alpha$, there is exactly one, with first component $S_m$, if $F(S_m)=\alpha$, and there are exactly two, with first components $\hat S$ and $S'$, if $F(S_m)<\alpha$.
		\end{enumerate}
	\end{enumerate}
	In all cases, the remaining components of the equilibria are given by \eqref{Q_i}.
\end{proposition}

\begin{proof}
	(i) In this case $D=(0,\infty)$ and, by \eqref{derivative_F_particular}, $F$ is strictly decreasing on $D$, with $F(S)\to0$ as $S\to\infty$. Hence, $F(S)=\alpha$ has a solution in $D$ if and only if $F(0)>\alpha$, and this solution is unique. In this case, the leading coefficient of $P$ is positive and $P(0)=m_1m_2(\alpha-F(0))<0$, so $P$ has a positive and a negative root, and the positive one is the largest, namely $\hat S$.
\medskip
    
	(ii) In this case $F$ is strictly increasing on $D$ by \eqref{derivative_F_particular}. Moreover, if $j$ is such that $S_{\max}=-m_j/A_j$, then $m_j+A_jS\to0^+$ as $S\to S_{\max}$, so that $F(S)\to\infty$. Hence, $F$ maps $D$ bijectively onto $(F(0),\infty)$, and $F(S)=\alpha$ has a solution $S^*\in D$ if and only if $F(0)<\alpha$, which is then unique. Since $F<\alpha$ on $[0,S^*)$ and $(m_1+A_1S)(m_2+A_2S)>0$ there, we have by \eqref{nova} that $P>0$ on $[0,S^*)$; in particular, $P(0)>0$. As the leading coefficient of $P$ is positive, $P$ is negative between its two roots. Therefore, $0$ cannot lie between the roots of $P$, and no root of $P$ lies in $[0,S^*)$, so that $S^*$ is the smallest root of $P$, namely $S'$.

\medskip
    
	(iii) Let $j\in\{1,2\}$ be such that $A_j<0$, so that $D=(0,-m_j/A_j)$ and $F(S)\to\infty$ as $S\to -m_j/A_j$. In case (a), the intermediate value theorem gives a solution of $F(S)=\alpha$ in $D$. On the other hand, the leading coefficient of $P$ is negative and $P(0)>0$, so that $P$ has exactly one positive root, which is the largest one, namely $S'$. Hence, $S'\in D$ is the unique solution in $D$.

	In cases (b) and (c), $F$ is strictly convex on $D$ by \eqref{derivative_F_particular}, so that $F'$ is strictly increasing on $D$. Moreover, $F'(S)\to\infty$ as $S\to -m_j/A_j$. In case (b), $F'>0$ on $D$ and therefore $F(S)>F(0)>\alpha$ for all $S\in D$. In case (c), $F'$ vanishes at a unique point $S_m\in D$, so that $F$ is strictly decreasing on $(0,S_m]$ and strictly increasing on $[S_m,-m_j/A_j)$. Since $F(0)>\alpha$ and $F(S)\to\infty$ as $S\to -m_j/A_j$, the equation $F(S)=\alpha$ has no solution in $D$ if $F(S_m)>\alpha$, the unique solution $S_m$ if $F(S_m)=\alpha$, and exactly two solutions if $F(S_m)<\alpha$. In the last case, both solutions are roots of $P$, so they are $\hat S$ and $S'$.
\end{proof}

Cases (i) and (ii) of \propref{prop:particular} correspond, for two phage types with the same latency period, to the situations described in \propref{Statement1}~(i) and (ii), respectively. In these cases, \propref{prop:particular} not only excludes coexistence equilibria for small or large doses, respectively, but also guarantees the existence of a unique biologically feasible coexistence equilibrium in the complementary case and provides its explicit expression. Case (iii), in which the two phage types behave differently, is not covered by \propref{Statement1}. In this case, if the doses are small, there is a unique biologically feasible coexistence equilibrium, while if the doses are large, there may be zero, one or two of them, coexisting with the locally asymptotically stable bacteria-free equilibrium.

\begin{rem}\label{rem:thresholds}
	The hypotheses on $A_1$ and $A_2$ in \propref{prop:particular} can be rewritten in terms of  the latency period $\tilde\tau$. Indeed, as in \remref{rem:n1}, each $A_i$ is an increasing function of $\tilde\tau$. If $b_i\leq1$, then $A_i>0$ for every $\tilde\tau>0$, while if $b_i>1$, then $A_i<0$ if and only if $\tilde\tau<\tilde\tau_i$, where
	\[
	\tilde\tau_i\coloneqq\frac{1}{\mu}\log\left(\frac{k_i^1b_i+\frac{k_i^2\alpha}{\mu}}{k_i^1+\frac{k_i^2\alpha}{\mu}}\right)\in\big(0,\tau_i^*\big).
	\]
	Setting $\tilde\tau_i\coloneqq0$ if $b_i\leq1$, we are in case (i) if $\tilde\tau>\max\{\tilde\tau_1,\tilde\tau_2\}$, in case (ii) if $\tilde\tau<\min\{\tilde\tau_1,\tilde\tau_2\}$, and in case (iii) if $\tilde\tau$ lies strictly between $\tilde\tau_1$ and $\tilde\tau_2$. In other words, long latency periods, for which coinfection makes both phage types unable to replace themselves, lead to coexistence only for large doses, short latency periods lead to coexistence only for small doses, and intermediate latency periods may lead to more complex situations.
\end{rem}


\section{Numerical Simulations}\label{Numerical Simulations}

In this section we illustrate the results of Section~\ref{Particular Case} by means of numerical simulations of system \eqref{system} with two phage types with the same latency period $\tilde\tau$. Some of the parameter values are taken from the estimations in \cite{BARDINA201399}, and the remaining ones have been chosen so that the conditions of the different cases in \propref{prop:particular} are satisfied. Time is measured in days, and the concentrations of bacteria and phages in units of $10^7$.

In all the simulations, the initial data are of the form $S_0(t)=c\,e^{t}$ for $t\in[-\tilde\tau,0]$, with $c>0$, and $Q_i^0(t)\equiv Q_i^0$ constant, while $I_0$ is given by \eqref{I0}. System \eqref{system} has been solved with a fourth-order Runge--Kutta scheme with constant step size $h=\tilde\tau/N$, $N\geq 20$, using linear interpolation for the delayed terms. In order to keep the accuracy when $S$ takes very small values, we integrate $\ln S$ instead of $S$. 

\subsection*{Case (i): the coexistence equilibrium as a threshold}

We take
\begin{eqnarray*}
&&(m_1,b_1,d_1,k_1^1,k_1^2,m_2,b_2,d_2,k_2^1,k_2^2,\alpha,\mu,\tilde\tau)\\&=&(0.18,3,0.4,27.36,25,0.1947,2.5,0.03,26,24,12.1622,0.5,0.1875).
\end{eqnarray*}
Then $A_1\approx7.045$ and $A_2\approx19.060$ are positive and $F(0)\approx64.806>\alpha$, so that, by \propref{prop:E0} and \propref{prop:particular}~(i), the bacteria-free equilibrium $E_0\approx(0,0,2.22222,0.15408)$ is locally asymptotically stable, with $\lambda_0=\alpha-F(0)\approx-52.644$, and there exists a unique biologically feasible coexistence equilibrium
\[
\hat E_c=(\hat S,\hat I,\hat Q_1,\hat Q_2)\approx(0.10599,\,0.23071,\,0.43165,\,0.01355).
\]
The characteristic equation at $\hat E_c$ has a real positive root $\lambda\approx2.916$, while all the other computed roots have negative real part, so that $\hat E_c$ is unstable.

\begin{figure}[ht]
	\centering
	\includegraphics[width=\textwidth]{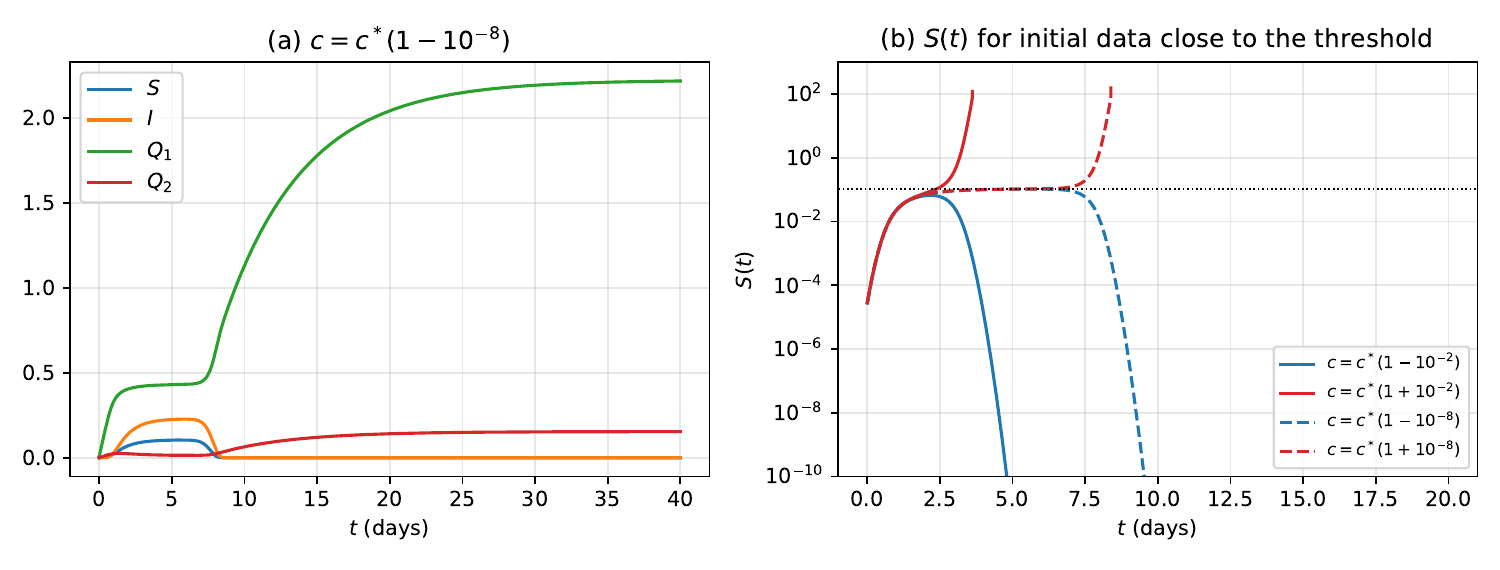}
	\caption{Case (i). (a) Trajectories of $S$, $I$, $Q_1$ and $Q_2$ for $c=c^*(1-10^{-8})$. (b) Evolution of $S$, in logarithmic scale, for initial data with $c=c^*(1\pm10^{-2})$ and $c=c^*(1\pm10^{-8})$; the dotted line corresponds to $\hat S$.}
	\label{fig:caseI}
\end{figure}

We take $Q_1^0=0.006$ and $Q_2^0=0.0005$. The simulations show a threshold value $c^*\approx2.8160\times10^{-5}$ such that, if $c<c^*$, the solution converges to the bacteria-free equilibrium, whereas if $c>c^*$ the phage concentrations tend to zero and $S$ grows exponentially without bound. Moreover, the closer $c$ is to $c^*$, the longer the solution remains close to $\hat E_c$ before leaving it (see Figure~\ref{fig:caseI}). That is, the unstable coexistence equilibrium separates the initial conditions leading to the eradication of the bacteria from those for which the therapy fails.

\subsection*{Case (ii): convergence to the coexistence equilibrium}

We take
\[
(m_1,b_1,d_1,k_1^1,k_1^2,m_2,b_2,d_2,k_2^1,k_2^2,\alpha,\mu,\tilde\tau)=(0.18,61,0.1,8,6,0.1947,50,0.15,5,4,12.1622,0.5,0.01875).
\]
Then $A_1\approx-474.08$ and $A_2\approx-241.76$ are negative and $F(0)\approx8.2965<\alpha$. Hence, the bacteria-free equilibrium $E_0\approx(0,0,0.55556,0.77042)$ is unstable, with $\lambda_0\approx3.866$, and, by \propref{prop:particular}~(ii), there exists a unique biologically feasible coexistence equilibrium
\[
E_c'=(S',I',Q_1',Q_2')\approx(1.52059\times10^{-4},\,3.45137\times10^{-5},\,0.92669,\,0.94974).
\]
The rightmost roots of the characteristic equation at $E_c'$ are $\lambda\approx-0.0533\pm0.8395\,i$, so that $E_c'$ is locally asymptotically stable and nearby solutions converge to it with damped oscillations of period approximately $2\pi/0.8395\approx7.5$ days.

\begin{figure}[ht]
	\centering
	\includegraphics[width=\textwidth]{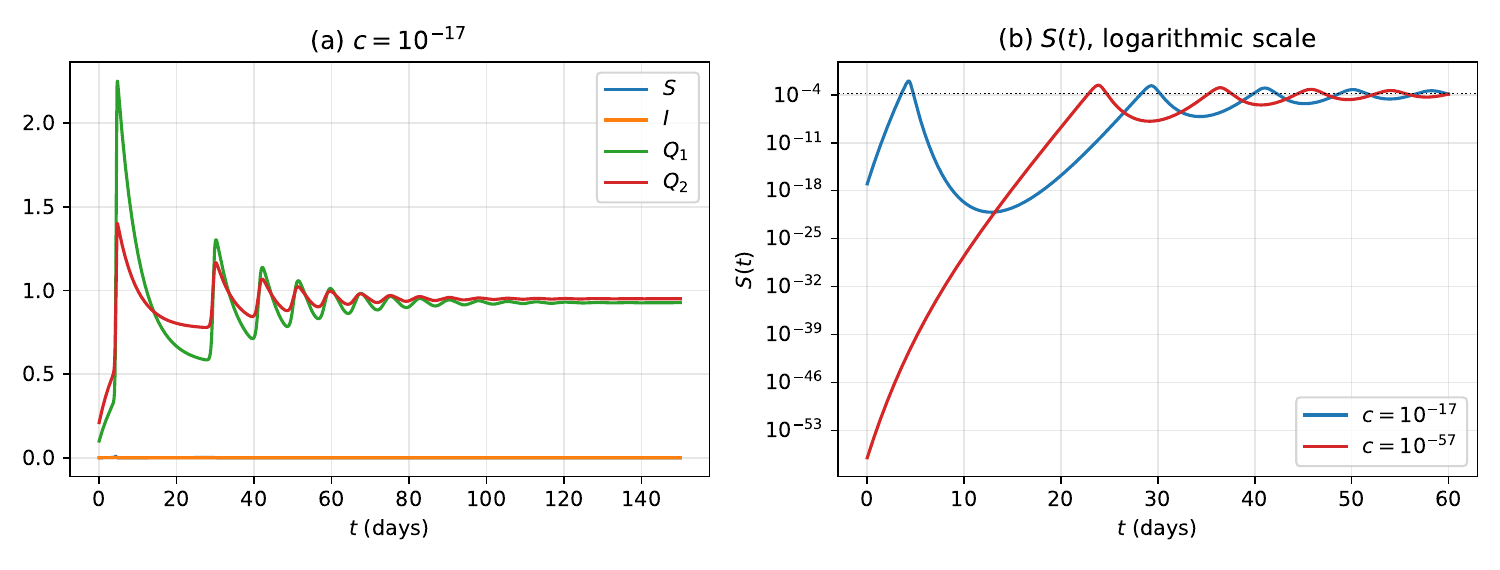}
	\caption{Case (ii). (a) Trajectories of $S$, $I$, $Q_1$ and $Q_2$ for $c=10^{-17}$. (b) Evolution of $S$, in logarithmic scale, for $c=10^{-17}$ and $c=10^{-57}$; the dotted line corresponds to $S'$.}
	\label{fig:caseII}
\end{figure}

Figure~\ref{fig:caseII}~(a) shows the solution with $Q_1^0=0.1$, $Q_2^0=0.21$ and $c=10^{-17}$, which converges to $E_c'$ with damped oscillations. Figure~\ref{fig:caseII}~(b) shows the evolution of $S$, in logarithmic scale, for $c=10^{-17}$ and $c=10^{-57}$. While $S$ is small, the phage concentrations evolve approximately as in the absence of bacteria, that is, $Q_i(t)\approx d_i/m_i+(Q_i^0-d_i/m_i)e^{-m_it}$, and therefore
\[
\ln S(t)\approx\ln c+\lambda_0t+\sum_{i=1}^2\frac{k_i^1}{m_i}\Big(\frac{d_i}{m_i}-Q_i^0\Big)\big(1-e^{-m_it}\big).
\]
Since $\lambda_0>0$, the solution eventually leaves the neighbourhood of $E_0$, and the time it remains close to $E_0$ grows like $|\ln c|/\lambda_0$ as $c\to0$, as can be observed in the figure. Notice also that, along the oscillations, $S$ takes extremely small values. For a real bacterial population, such values would correspond to the extinction of the bacteria, which illustrates a well-known limitation of deterministic models.

\subsection*{Case (iii)(c): bistability}

Finally, we take
\begin{eqnarray*}
&&(m_1,b_1,d_1,k_1^1,k_1^2,m_2,b_2,d_2,k_2^1,k_2^2,\alpha,\mu,\tilde\tau)\\&=&(0.18,13.7,0.39,1.8,19.4,0.1947,2.77,0.67,4.15,19.5,12.1622,0.5,0.093).
\end{eqnarray*}
Then $A_1\approx-0.2989<0<A_2\approx14.728$, $F(0)\approx18.181>\alpha$ and $F'(0)\approx-1073.8<0$, and the minimum of $F$ on $D=(0,-m_1/A_1)\approx(0,0.60225)$ is $F(S_m)\approx6.287<\alpha$, attained at $S_m\approx0.1227$. By \propref{prop:particular}~(iii)(c), there exist exactly two biologically feasible coexistence equilibria,
\[
\hat E_c\approx(0.00981,\,0.01084,\,2.20254,\,1.97533)
\qquad\text{and}\qquad
E_c'\approx(0.40162,\,0.44387,\,6.50395,\,0.10966).
\]
Notice that both phage types are effective, since $b_1e^{-\mu\tilde\tau}\approx13.08$ and $b_2e^{-\mu\tilde\tau}\approx2.64$, but $A_2>0$ and that implies (see \remref{rem:R0}, and notice that $A_2=B_2$) that $R_{0,2}<1$. That is, coinfection prevents phages of type $2$ from replacing themselves at the coexistence equilibria. The bacteria-free equilibrium $E_0\approx(0,0,2.16667,3.44119)$ is locally asymptotically stable, with $\lambda_0\approx-6.019$. The characteristic equation at $\hat E_c$ has a real positive root $\lambda\approx0.858$, so that $\hat E_c$ is unstable, while the rightmost roots at $E_c'$ are $\lambda\approx-0.440\pm0.649\,i$, so that $E_c'$ is locally asymptotically stable.

\begin{figure}[ht]
	\centering
	\includegraphics[width=\textwidth]{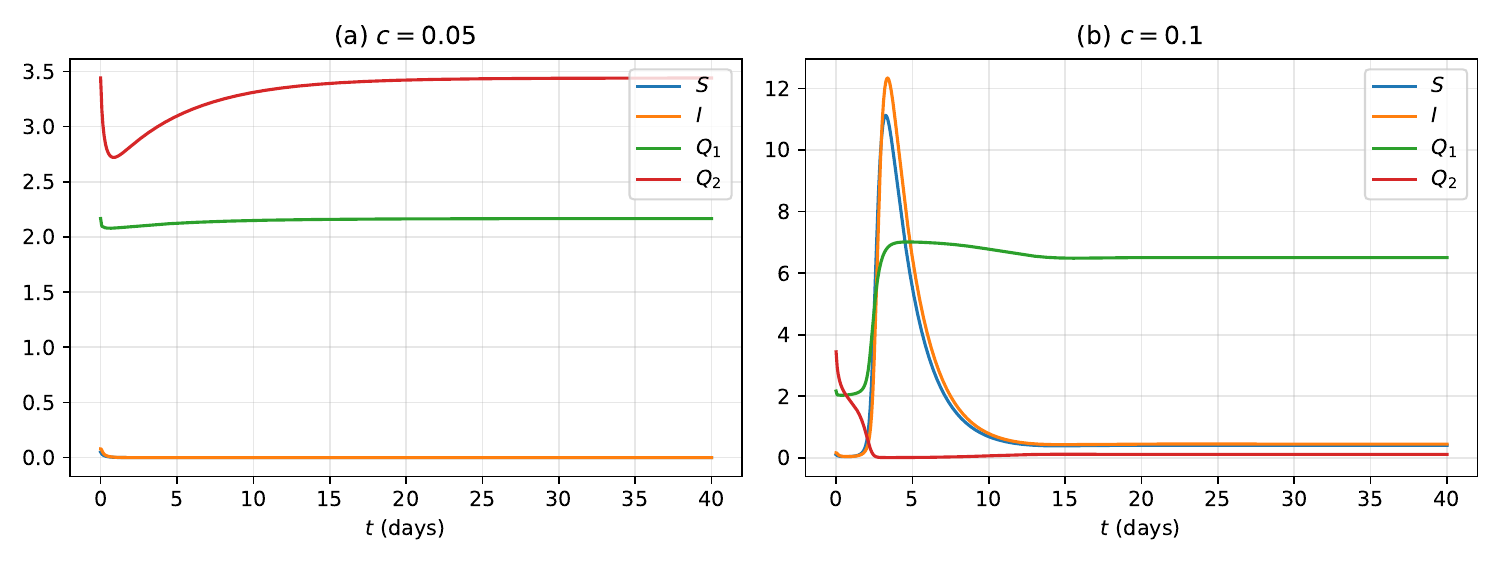}
	\caption{Case (iii)(c). Trajectories of $S$, $I$, $Q_1$ and $Q_2$ for (a) $c=0.05$ and (b) $c=0.1$.}
	\label{fig:caseIIIc}
\end{figure}

We take $Q_i^0=d_i/m_i$, that is, the phage concentrations at the moment the bacteria are introduced are those of the bacteria-free equilibrium. Figure~\ref{fig:caseIIIc} shows that, for $c=0.05$, the solution converges to the bacteria-free equilibrium, whereas for $c=0.1$ it converges to the coexistence equilibrium $E_c'$. Hence, in this case the outcome of the therapy depends on the initial bacterial load: small infections are eradicated, while larger ones persist at a positive level, even though the administered doses are large.

\section{Concluding Remarks}\label{Conclusions}

In this paper we have analyzed a model for phage therapy in which a cocktail of $n$ phage types, each with its own latency period, is administered to a bacterial population, and in which coinfection is allowed. This model extends the single-phage model with coinfection in \cite{bardina2019coinfection} and the coinfection-free cocktail model in \cite{bascompte2014models}.

We have shown that the initial value problem has a unique global nonnegative solution and that, under suitable conditions on the parameters and the initial data, there exists an invariant region in which the susceptible bacteria are nonincreasing and, under a slightly stronger condition, are eradicated exponentially fast. This invariant region only exists when the bacteria-free equilibrium is locally asymptotically stable, which happens if and only if the doses are large, in the sense that $\sum_{i=1}^n k_i^1d_i/m_i>\alpha$.

Regarding coexistence equilibria, we have obtained sufficient conditions for their nonexistence, which can be interpreted in terms of the expected number $R_{0,i}$ of phages of type $i$ produced per adsorbed phage of type $i$ at a coexistence equilibrium. If $R_{0,i}<1$ for every phage type, coexistence equilibria can only exist when the doses are large, whereas if $R_{0,i}>1$ for every phage type, they can only exist when the doses are small. For a single phage type, this analysis is complete, and shows that coinfection lowers the latency period above which coexistence occurs for large doses. The main consequence of coinfection is that phages lost by adsorbing to infected bacteria reduce $R_{0,i}$, so that effective phages may behave as ineffective ones. As a consequence, in contrast with the coinfection-free case, coexistence equilibria may exist even when the doses are large.

For a cocktail of two phage types with the same latency period, we have obtained a complete description of the biologically feasible coexistence equilibria, together with explicit expressions. In this case, the quantity $R_{0,i}$ only depends on the parameters of the model. This number takes into account both the probability that the phage adsorbs to a susceptible bacterium rather than to an infected one, in which case it is lost, and the probability that the infected bacterium survives the latency period, but not the natural decay of phages. Thus, $R_{0,i}<1$ means that phages of type $i$ cannot replace themselves through infection, while $R_{0,i}>1$ means that they can. If $R_{0,1}<1$ and $R_{0,2}<1$, there exists a biologically feasible coexistence equilibrium if and only if the doses are large, and it is unique. If $R_{0,1}>1$ and $R_{0,2}>1$, there exists such an equilibrium if and only if the doses are small, and it is unique. Finally, if one of the two quantities is smaller than one and the other one is larger than one, there is a unique such equilibrium when the doses are small, while when the doses are large there may be zero, one or two of them.

The numerical simulations illustrate these results and suggest that the outcome of the therapy may depend on the initial conditions. In the first example, the unique coexistence equilibrium is unstable and acts as a threshold: depending on the initial bacterial concentration, the bacteria are either eradicated or grow without bound. In the second example, the bacteria-free equilibrium is unstable and the solutions converge to a stable coexistence equilibrium with damped oscillations, although the bacterial concentration reaches values that, for a real population, would correspond to extinction. In the third example, with two coexistence equilibria, the system is bistable: small initial bacterial loads are eradicated, while larger ones persist at a positive level, even though the doses are large and both phage types are effective.

A natural continuation of this work is the study of a stochastic version of the model. As in \cite{BARDINA201399,bardina2019coinfection}, several random effects, such as the noise in the data collected from laboratory experiments, random fluctuations of the parameters due, for instance, to changes in temperature, or the randomness of the latency periods, can be summarized in a small random perturbation of system \eqref{system} driven by a Wiener process. For the single-phage models in \cite{BARDINA201399,bardina2019coinfection}, the exponential convergence of the deterministic system to the bacteria-free equilibrium is the key ingredient to prove that, in a small noise regime, the solutions of the perturbed system concentrate around this equilibrium.  Therefore, it would be interesting to obtain a similar concentration result for a stochastic version of the cocktail model, and to study how the composition of the cocktail and coinfection affect the convergence to the bacteria-free equilibrium.

\printbibliography
\end{document}